\documentclass[12pt,reqno]{amsart} 
\usepackage{mathrsfs}
\usepackage{cite}
\usepackage{amssymb,amsmath}
\def\re#1{\par\hang\tex{#1}}

\def\cal{\mathcal}
\def\gg{\frak g}

\def\d{\delta}
\def\D{\Delta}
\def\UU{{U}}

\def\e{\epsilon}

\def\vp{\varepsilon}

\def\es{\varepsilon}

\def\v{\varphi}

\def\sc{\scriptstyle}
\def\ssc{\scriptscriptstyle}

\def\cl{\centerline}

\def\rar{\rightarrow}

\def\vs{\vspace*}

\def\ni{\noindent}
\def\ptl{\partial}

\def\N{\mathbb{N}{\ssc\,}}
\def\Z{\mathbb{Z}{\ssc\,}}

\def\C{\mathbb{C}{\ssc\,}}
\def\F{\mathbb{F}{\ssc\,}}
\def\QED{\hfill$\Box$} \numberwithin{equation}{section}
\newtheorem{theo}{Theorem}[section]
\newtheorem{conv}[theo]{Convention}
\newtheorem{defi}[theo]{Definition}
\newtheorem{rema}[theo]{Remark}
\newtheorem{lemm}[theo]{Lemma}
\newtheorem{prop}[theo]{Proposition}

\def\adddot{$\!\!\!${\bf.}\ \ }
\def\adddot{}

\begin{document}

\cl{{\large \bf Dual Lie Bialgebras of  Witt and Virasoro Types}\footnote{
Supported by NSF grant  11071147, 10825101
 of China,
NSF grant ZR2010AM003 of Shandong Province,
SMSTC grant 12XD1405000, and Fundamental Research Funds for the Central Universities.
}}
\vs{6pt}

\cl{ Guang'ai Song$^{\,*}$, Yucai Su$^{\,\dag}$}

\cl{\small $^{*\,}$College of Mathematics and Information Science, Shandong
Institute of \vs{-4pt}Business }

\cl{\small and Technology, Yantai, Shandong 264005, China}

\cl{$^{\dag\,}$Department of Mathematics, Tongji University, Shanghai 200092, China}

\cl{\small E-mail: gasong@sdibt.edu.cn, ycsu@tongji.edu.cn}\vskip5pt

\noindent{\small{\bf Abstract.} Structures of dual Lie bialgebras on the one sided
 Witt algebra, the Witt algebra and the Virasoro algebra are
investigated. As a result, we obtain some infinite dimensional Lie algebras.


}

\section{Introduction}\setcounter{section}{1}\setcounter{equation}{0}

 Lie bialgebras,  subjects of
intensive study in recent literature (e.g., \cite{1,1-2,D,G,M1,3,4,5,6,7,8,9,10,12,13,14}), 
are important ingredients in
quantum groups \cite{2}, which have close relations with
the Yang-Baxter equation.   In general, a Lie
bialgebra is a vector space endowed with both the structure of Lie
algebra and the structure of Lie coalgebra simultaneously,
satisfying some compatibility condition, which was suggested by a
study of Hamiltonian mechanics and Poisson Lie group \cite{2,3}.
Michaelis \cite{4} studied a class of Witt type Lie bialgebras, and  introduced   techniques  on how to construct
 coboundary or triangular Lie bialgebras from Lie algebras which
 contains two linear independent elements $a$ and $b$ satisfying condition
 $[a, b] = kb$ for some non-zero scalar $k$. Ng and Taft \cite{9} proved that all structures of Lie bialgebras  on the one sided
 Witt algebra, the Witt algebra and the Virasoro algebra are
 coboundary triangular (see also \cite{8}). Furthermore, they gave a complete classification of all structures of
 Lie bialgebras  on the one sided
 Witt algebra.
  For the cases of generalized Witt type Lie algebras, the authors \cite{10} obtained that
 all  structures of Lie bialgebras on them are
 coboundary triangular. Similar results also hold
 for some other kinds of Lie algebras (see, e.g., \cite{13, 14}).

It may sound that
coboundary triangular Lie bialgebras have relatively simple structures.
However, even for the (two-sided) Witt algebra and the Virasoro algebra,
a classification of coboundary triangular Lie bialgebra structures on them is still an open problem.
Thus it seems to us that it is worth
paying more attention on them.
This is one of our motivations in the present paper, here we investigate Lie bialgebra structures from the point of view of dual Lie bialgebras.
Investigating dual structures on some algebraic structures has also drawn some authors' much attention, e.g.,
\cite{1-2,D,G,M1,5,6}.
One may have noticed that the dual of a finite dimensional Lie bialgebra is naturally a Lie bialgebra, and so one would not predict anything new in this case.
Thus we start our investigation by considering the ``dual'' of  the one sided
 Witt algebra, the Witt algebra and the Virasoro algebra. Notice that the (full) dual of an infinite dimensional space is always an infinite dimensional space with uncountable
 basis, thus here  ``dual'' means some restricted dual (cf.~\eqref{max-good}).
A surprising thing is that dualizing Lie bialgebras may produce new Lie algebras, as can be seen from Remark \ref{final-rem}. This is another motivation in the present paper.

The paper is organized as follows. Some
 definitions and preliminary results are briefly recalled in Section 2. Then in Section 3, structures of dual coalgebras
 of $\F[x]$ and $\F[x^{\pm 1}]$ are addressed. Finally in Section 4,  structures of
 dual Lie bialgebras of the one side Witt Lie bialgebras, the Witt Lie
 bialgebras and the Virasoro Lie bialgebras are investigated. As a result, we obtain some infinite dimensional Lie algebras.
 The main results of the present paper are summarized in Theorems \ref{Theo3.1}, \ref{dual-L} and
 \ref{Last-Theo}.
 \vskip5pt

\section{Definitions and preliminary
results}\setcounter{section}{2}\setcounter{equation}{0}

Throughout the paper, $\F$ denoted an algebraically closed field of characteristic zero unless otherwise stated.
All vector spaces are assumed to be over $\F$. As usual,
we use $\Z$, $\Z_+$, $\N$ to denote the sets of integers, nonnegative
integers, positive integers, respectively.

We briefly recall some notions on Lie bialgebras, for details, we refer readers to, e.g., \cite{2,10}.
\begin{defi}\adddot{\label{2.1}}\rm\begin{enumerate}\parskip-1pt\item A {\it Lie bialgebra}
is a triple $(L, [\cdot, \cdot], \d)$ such that\begin{enumerate}
\item[(i)]  $(L, [\cdot, \cdot])$ is  a   Lie  algebra;
\item[(ii)] $(L, \d)$ is   a  coalgebra;
\item[(iii)] $\d[x, y] = x \cdot \d(y) - y\cdot\d(x)$ for  $x,y \in L$,
\end{enumerate}
\noindent where  $ x \cdot (y\otimes z) = [x, y]\otimes z + y
\otimes [x, z]$ for  $x, y, z \in L.$
\item A  Lie bialgebra
 $(L, [\cdot, \cdot],  \d)$ is {\it coboundary} if $\d$ is coboundary in the sense that there exists $r \in L \otimes L$ written as $r= \sum r^{[1]} \otimes r^{[2]} $, such that $\d(x) = x\cdot r$
for $x\in L$.
\item A coboundary Lie bialgebra
 $(L, [\cdot, \cdot],  \d)$ is {\it triangular} if $r$ satisfies  the following {\it classical Yang-Baxter Equation} (CYBE),
\begin{equation}\label{CYBE}C(r) = [r_{12}, r_{13}] + [r_{12}, r_{23}] + [r_{13}, r_{23}]
 =0,\end{equation}
 \noindent where
$r_{12} = \sum r^{[1]} \otimes r^{[2]} \otimes 1,$
$r_{13}= \sum r^{[1]} \otimes 1 \otimes r^{[2]},$
$r_{23} = \sum r^{[1]} \otimes 1 \otimes r^{[2]}$
are elements in $ \UU(L)\!
\otimes\! \UU(L) \!\otimes\! \UU(L)$, and $\UU(L)$ is the universal enveloping algebra of $L$.\end{enumerate}
 \end{defi}

Two Lie bialgebras $(\gg,[\cdot,\cdot],\d)$ and $(\gg',[\cdot,\cdot]',\d')$ are said to be {\it dually paired}
if there bialgebra structures are related via
\begin{equation}\label{Dual-paired}
\langle[f, h ]', \xi\rangle = \langle f \otimes h, \d \xi \rangle ,\ \ \  \langle \d' f, \xi \otimes
\eta \rangle = \langle f, [\xi, \eta]\rangle
\mbox{ \ for $f, h \in \gg',\ \xi, \eta \in \gg,$ }\end{equation}
 where $\langle\cdot,\cdot\rangle$ is a nondegenerate
bilinear form on $\gg'\times\gg$, which is naturally extended to a nondegenerate
bilinear form on \mbox{$(\gg'\otimes\gg')\times(\gg\otimes\gg)$.}
In particular, if $\gg'=\gg$ as a vector space, then $\gg$ is called a {\it
self-dual Lie bialgebra}.

The following easily obtained result can be found in \cite{3}.
\begin{prop}\adddot
Let $(\gg, [\cdot, \cdot], \d)$ be a finite dimensional Lie
bialgebra, then so is the linear dual space $\gg ^*:={\rm Hom}_\F(\gg ,\F)$ by dualisation, namely $(\gg^*,[\cdot,\cdot]',\d')$ is the Lie bialgebra
defined by  \eqref{Dual-paired} with $\gg'=\gg^*$. In particular, $\gg$ and $\gg^*$ are dually paired.
\end{prop}

Thus a finite dimensional $(\gg,[\cdot,\cdot],\d)$ is always self-dual as there exists a vector space isomorphism $\gg\to\gg^*$ which pulls back the bialgebra structure on $\gg^*$ to $\gg$ to obtain another bialgebra structure on $\gg$ to make it to be self-dual.
However, in sharp contrast to finite dimensional case, infinite
dimensional Lie bialgebras are  not self-dual in general.

For convenience, we denote by $\v$  the Lie bracket of Lie algebra $(\gg ,
[\cdot, \cdot])$, which can be regarded as a linear map $\v:\gg\otimes\gg\to\gg$. Let $\v^{\ast}:\gg^*\to(\gg\otimes\gg)^*$ be  the dual of $\v$.

\begin{defi}\adddot\rm \cite{M1}
Let $(\gg , \v)$ be a Lie algebra over $\F$. A  subspace $V$ of $\gg ^*$ is called a {\it good
subspace} if $\v^*(V) \subset V \otimes V.$
Denote $\Re = \{ V\, |\, V {\rm\ is \ a \ good \ subspace \ of \ }\gg ^*\}$. Then
\begin{equation}\label{max-good}\gg ^{\circ} = \mbox{$\sum\limits_{V\in \Re}$} V,\end{equation}
 is also a good subspace
of $\gg ^*$, which is obviously the maximal good subspace
of $\gg ^*$.
\end{defi}

\begin{prop}\adddot{\rm\cite{M1}}\label{PPPP} For any good subspace $V$ of $\gg ^*$, the pair $(V,
\v^{\ast})$ is a Lie coalgebra. In particular, $(\gg ^{\circ}, \v^*)$
is a Lie coalgebra.
\end{prop}
It is clear that
if $\gg $ is a finite dimensional Lie algebra, then
$\gg ^{\circ} = \gg ^*$. However if $\gg $ is infinite dimensional, then
$\gg ^{\circ} \subset \gg ^* $ in general.

 For any Lie algebra $\gg $, the dual space $\gg ^*$ has a natural right $\gg $-module
 structure defined for $f\in \gg ^*$ and $ x\in \gg $ by $$(f\cdot x)(y) = f([x,
 y])\mbox{ \ for \ }y\in \gg .$$ We denote $f\cdot \gg  = {\rm span} \{ f\cdot x\, |\, x \in
 \gg \}$, the {\it space of translates} of $f$ by elements of $\gg $.

We summarize some results of \cite{1,1-2,D} as follows.
 \begin{prop}\adddot\label{Th111}
 Let $\gg $ be a Lie algebra. 
 Then\begin{enumerate}\item $\gg ^{\circ}=\{f\in\gg^*\,|\,f\cdot\gg\mbox{ is finite dimensional}{\sc\,}\}.$
\item
$\gg ^{\circ}=(\v^*)^{-1}(\gg^*\otimes\gg^*)$, the preimage of $\gg^*\otimes\gg^*$ in $\gg^*$.
\end{enumerate}
 \end{prop}

The notion of good subspaces of an associative algebra can be defined analogously.
In the next two sections, we shall investigate $\gg^\circ$ for some associative or
Lie algebras $\gg$.\vskip10pt

\section{The structure of $\F[x,x^{-1}]^{\circ}$
}\setcounter{section}{3}\setcounter{theo}{0}\setcounter{equation}{0}

Let $(\cal{A}, \mu, \eta)$ be an associative $\F$-algebra with unit,  where $\mu$ and $\eta$ are respectively the multiplication $\mu: \cal{A} \otimes
\cal{A} \rar \cal{A}$ and  the unit $\eta: \F \rar \cal{A}$,
satisfying $$\begin{array}{ll}
\mu \circ (id \otimes \mu)= \mu \circ (\mu \otimes id)):&
\cal{A} \otimes \cal{A} \otimes \cal{A} \rar \cal{A},\\[2pt]
(\eta\otimes id)(k\otimes a) = (id \otimes \eta)(a \otimes k): &\F \otimes
\cal{A}\cong \cal{A}\otimes \F \cong \cal{A},\end{array}$$
 for  $k \in \F,\, a, b, c \in \cal{A}$.
Then a coassociative coalgebra is a triple $(C, \D, \e)$,  which is obtained
by conversing arrows in the definition of an associative algebra.
Namely,  $\D: C \rar C\otimes C$ and $\e: \F \rar C$ are respectively
comultiplication and counit of $C$, satisfying
$$\begin{array}{ll}(\D \otimes id)\circ \D = (id \otimes \D)\circ \D: &C \rar
C\otimes C \otimes C,\\[2pt]
(\e \otimes id)\circ \D = (id \otimes
\e)\circ \D: &C \rar C\otimes C \rar \F\otimes C \cong C\otimes \F
\cong C.\end{array}$$

For any vector space $\cal{A}$, there exists a natural  injection $\rho:
\cal{A^{\ast}} \otimes \cal{A^{\ast}} \rar (\cal{A} \otimes
\cal{A})^{\ast}$ defined by $\rho(f, g )(a, b)= \langle f, a\rangle \langle g, b\rangle $ for
$f, g \in \cal{A^{ \ast}}$ and $ a, b, \in \cal{A}$. In case $\cal A$ is finite dimensional,
$\rho$ is an isomorphism.
For any
algebra $(\cal{A}, \mu),$ the multiplication $\mu: \cal{A} \otimes
\cal{A} \rar \cal{A}$ naturally induces the map $\mu^{\ast}: \cal{A^{\ast}}
\rar (\cal{A} \otimes \cal{A})^{\ast}$. If $\cal{A}$ is  finite
dimensional, then the isomorphism
$\rho$ insures that $(\cal{A^{\ast}}, \mu^{ \ast},
\e^{\ast})$ is a coalgebra, where for simplicity, $\mu^{\ast}$ denotes the
composition of the maps: $ \cal{A^{\ast}}
\stackrel{\mu^{\ast}}{\rar}(\cal{A} \otimes \cal{A})^{\ast}
\stackrel{(\rho)^{{\rm -1}}}{\rar} \cal{A^{\ast}} \otimes
\cal{A^{\ast}}$. However, if $\cal{A}$ is infinite
dimensional, the situation is different. 

Let $\cal A=\F[x]$ (resp., $\F[x,x^{-1}]$) be the algebra of polynomials (resp., Laurent polynomials) on variable $x$.
As a vector space, we have ${\cal A}^*=\F[[\vp]]$ (resp., $\F[[\vp,\vp^{-1}]]$), the space of formal power series on $\vp$ (resp., $\vp,\vp^{-1}$), where $\vp^i$ is the dual element of $x^i$, namely, $\langle\vp^i,x^j\rangle:=\vp^i(x^j)=\d_{i,j}$ for $i,j\in\Z_+$ (resp., $i,j\in\Z$). In this way, for any $f=\sum_if_i\vp^i\in\cal A^*$ (which can be an infinite sum) and any $g=\sum_jg_jx^j\in{\cal A}$ (a finite sum), we have
\begin{equation}\label{Dualll}
\mbox{$\langle f,g\rangle=f(g)=
\sum\limits_{i,j}f_i\langle\vp^i,g\rangle
=\sum\limits_{i,j}f_ig_j\langle\vp^i,x^j\rangle=\sum\limits_{i,j}f_ig_j\vp^i(x^j)=\sum\limits_j f_jg_j,
$}
\end{equation}
where the right-hand side is a finite sum.

We summarize some results of \cite{7} as follows.

\begin{prop}\adddot\label{F[x]}{\rm\cite{7}}
For any $f= \sum_{i=0}^{\infty} f_i \vp^i \in
 \F[[\vp]]$ with $f_i\in\F$, \begin{equation}\label{poly}f \in \F[x]^{\circ}\ \Longleftrightarrow\ \exists\,r\in\N \mbox{ \ such that }f_n \!=\! h_1f_{n-1}\! +\! h_2
f_{n-2}
\! +\! \cdots
\! +\! h_r f_{n-r},\end{equation}
 for some $h_i\in\F$ and all $ n > r$.
\end{prop}

One can easily generalize this result to the Laurent polynomial algebra $\F[x,x^{-1}]$ as follows.

\begin{theo}\adddot \label{F[x^{pm1}]}
For any $f = \sum_{j=-\infty}^{\infty}f_j\es^j \in
 \F[[\vp,\vp^{-1}]]$ with $f_i\in\F$,
 \begin{equation}\label{lau-poly}f\in\F[x,x^{-1}]^\circ\ \Longleftrightarrow\ \exists\,r\in\N\mbox{ \ such that }f_n \!=\! h_1f_{n-1}\! +\! h_2 f_{n-2}
 \!+\! \cdots 
 \!+\! h_r f_{n-r},\end{equation}
for some $h_i\in\F$ and all $n \in \Z$.
\end{theo}

\begin{rema}\adddot\rm\label{Poly-and-Laupoly}
Note that there is a significant difference between $\F[x]^\circ$ and $\F[x,x^{-1}]^\circ$. From \eqref{poly}, we can easily see that
every polynomial of $\vp$ is in $\F[x]^\circ$, namely, $\F[\vp]\subset\F[x]^\circ$. In fact one can prove that  $\F[x]^\circ$ is the subspace of $\F[[\vp]]$ consisting of rational functions $r(\vp)$ on $\vp$ such that zero is not a singular point of $r(\vp)$ , i.e., $\F[x]^\circ=\{\frac{f(\vp)}{g(\vp)}\,|\,f(\vp),g(\vp)\in\F[\vp],\,g(0)\ne0\}$. However, from \eqref{lau-poly},
one can obtain that if $0\ne f\in\F[x,x^{-1}]^\circ$, then for any $N\in\N$, there exist $i>N$ and $j<-N$ such that $f_i\ne0\ne f_j$. In particular, $\F[\vp,\vp^{-1}]\cap\F[x,x^{-1}]^\circ=0$.
\end{rema}

\ni{\it Proof of Theorem \ref{F[x^{pm1}]}.~}~Note that $\F[x,x^{-1}]$ is a principle ideal domain.
Let $f \in \F[x^{\pm 1}]^{\circ}$. By \cite{12},  ${\rm Ker}(f)
\supset I$, where $I$ is a finite codimensional ideal of $\F[x,x^{-
1}]$. Thus there exists a polynomial $h(x) = x^r - h_1 x^{r-1} - \cdots-h_r \in
\F[x]$ with $ r\in\N$ and $h_i\in\F$ such that $I = (h(x)) $. Then
for $n \in \Z $, we have $ x^{n-r}h(x) \in I \subseteq {\rm Ker}(f), $ which means
$f(x^{n-r}h(x)) =0,$ i.e., $f_n = h_1f_{n-1} + h_2 f_{n-2}
 + \cdots + h_r f_{n-r} $.

Now assume  $f_n = h_1f_{n-1} + h_2 f_{n-2}
 + \cdots  + h_r f_{n-r} $ for  $n \in \Z$. Set $h(x) = x^r - h_1x^{r-1} - \cdots - h_r$.
Then
$ f(x^{n-r}h(x)) = 0$, i.e., the ideal $I=(h(x))$ is a finite codimensional ideal such that $I\subset{\rm Ker}(f)$. Thus, $f \in
\F[x,x^{- 1}]^{\circ}$ by Proposition \ref{Th111}.
 \QED\vskip5pt

A sequence $\{f_n \}_{n=-\infty}^{\infty}$
satisfying the condition $f_n = h_1f_{n-1} + h_2 f_{n-2}
 + \cdots + h_r f_{n-r}$ for $n \in \Z$ is called a {\it recursive
 sequence}, the minimal integer $r$ of the recursive relation is called the
{\it degree} of the sequence. For convenience, we call a formal power series $f=\sum_{n=-\infty}^{\infty} f_n \es^n$
a {\it  recursive power series} whenever $\{f_n \}_{n=-\infty}^{\infty}$ is a
 recursive sequence.

A coalgebra $C$ over a field $\F$ is called
{\it irreducible} if any two nonzero subcoalgebras of $C$ have nonzero
intersection. A subcoalgebra of $C$ is called an {\it irreducible
component} of $C$ if it is a maximal irreducible subcoalgebra.

The following result can be found in \cite{12}.
\begin{prop}\adddot If $C$ be a cocommutative coalgebra, then $C$ is direct sum of its irreducible component.
\end{prop}

\def\C{\F}
Let $\cal{A}=\C[x,x^{-1}]$ (or $\C[x])$.
 For $a\in\F$, let  $I_a$ be the ideal generated by $(x-a)$ in
$\cal{A}$. Then
$U = \{I_a\,|\, 0\ne a\in \C\}$ (or
$U = \{I_a \,|\, a\in \C\}$) is the set of maximal ideals of
$\cal{A}$. By \cite{12}, $I_a^{\perp} := \{f \in \cal{A}^{\ast}\, |\,f(I_a)=0\}$ is a simple subcoalgebra of $\cal{A^{\circ}}$. For $n\in\N$, denote
$I_{a^n} = ((x-a)^n).$ Then we have the chain of ideals of $\cal{A}$:  $I_a
\supset I_{a^2} \supset \cdots \supset I_{a^n} \supset \cdots. $
Since each ideal $I_{a^n}$ is a finite codimensional ideal of $\cal{A}$, we have
$$I_a^{\perp} \subset I_{a^2}^{\perp} \subset \cdots \subset
I_{a^n}^{\perp} \subset \cdots,$$ which is a chain of subcoalgebras of
$\cal{A^{\circ}}$. We denote $S_a = \cup_{n=0}^{\infty}
I_{a^n}^{\perp}= \lim_{n\rar \infty}I_{a^n}^{\perp}$. Then
$S_a$ is a maximal irreducible subcoalgebra of $\cal{A^{\circ}}$.
\begin{theo}\adddot\label{Theo3.1} Let $S_a$ be defined as
above. Then $\cal{A^{\circ}} $$=
\oplus_{a \in \C^{*}} S_a$ if $\cal{A}= \C[x,x^{-1}]$, or
$\cal{A^{\circ}}= \oplus_{a \in \C} S_a$ if $\cal{A}= \C[x]$.
\end{theo}

\ni{\it Proof.~}~We prove the case $\cal{A} = \C[x,x^{-1}]$ as  the
other case is analogous. For  $0 \neq a \in \C$, $S_a$ is
obviously a maximal irreducible subcoalgebra of $\cal{A^{\circ}}$.
If $S_a \cap S_b \neq 0 $ for some $a,b\in\F$ with $a\ne b$,
 it must contain a
simple subcoalgebra of $\cal{A^{\circ}}$. Since $I_a^{\perp}$ is a
unique simple subcoalgebra of $S_a$, we obtain $I_a^{\perp}\subset S_a\cap S_b$. Analogously,
$I_b^{\perp}\subset S_a\cap S_b$, which forces  $I_a^{\perp}=
I_b^{\perp},$ contradiction with $a\neq b$.
Thus the sum is directed.
Since the irreducible
components of $\cal{A^{\circ}}$ are in  1-1 correspondence with the finite codimensional maximal ideals of
$\cal{A}$, we see that  $S_a$'s with $0\neq a \in \C$ run over all
the irreducible components of $\cal{A^{\circ}}$.
\QED

\section{Dual Lie bialgebras of Witt type
}\setcounter{section}{4}\setcounter{theo}{0}\setcounter{equation}{0}

Let $(\gg , \v, \d)$  be a Lie bialgebra, and $\v^*: \gg ^* \rar
(\gg \otimes \gg )^*$ the dual of $\v$. By Proposition \ref{PPPP},
$\v^*$ induces a map $\v^{\circ}:=\v^*|_{\gg^{\circ}}:\gg ^{\circ}
\rar \gg ^{\circ} \otimes \gg ^{\circ}$, which makes
$(\gg^{\circ},\v^\circ)$  to be a Lie coalgebra. By \cite[Proposition
3]{8}, the map $\d^*:\gg ^* \otimes \gg ^* \hookrightarrow (\gg
\otimes \gg )^* \stackrel{\d^*}{\rar} \gg ^*$ induces a map
 $\d^{\circ}:=\d^*_{\gg ^{\circ} \otimes \gg ^{\circ}}:
\gg ^{\circ} \otimes \gg ^{\circ} \rar \gg ^{\circ}$, which makes
$(\gg^{\circ},\d^\circ)$ to be a Lie algebra. Thus we obtain a Lie bialgebra
  $(\gg ^{\circ}, \d^{\circ}, \v^{\circ} )$, called  the {\it dual Lie
 bialgebra} of $(\gg , \v,  \d)$.

\begin{conv}\rm When there is no confusion, we will use $[\cdot,\cdot]$ to denote the bracket in $\gg$ or  $\gg^\circ$, i.e., $[\cdot,\cdot]=\v$ or $\d^{\circ}$, and we also use $\D$ to denote the cobracket in $\gg$ or  $\gg^\circ$, i.e., $\D=\d$ or
 $\v^\circ$.
\end{conv}
If $\gg $ is an infinite dimensional Lie algebra, there do not exist effective
approaches to describe the structure of $\gg ^{\circ}$ explicitly. However,
if $\gg $ is a Lie algebra ``induced'' from some associative algebra (in sense of \eqref{AAAA} below), it is sometimes possible to
use the theory of dual associative algebras to describe $\gg ^{\circ}$.
So, let $(\cal{A}, \mu)$ be a commutative associative $\F$-algebra with multiplication $\mu$.
Let $\partial\in {\rm Der}_{\F} (\cal{A})$ be a derivation of $\cal A$.
Then we obtain a Lie algebra, denoted ${\cal A}_\ptl$, on the space $\cal A$, with bracket defined by
\begin{equation}\label{AAAA}id \otimes \ptl - \ptl\otimes id,\mbox{ \ namely, \ }[a, b] = a\ptl(b) -
\ptl(a)b\mbox{ \ for \ }a,b\in{\cal A}.\end{equation}
Such a Lie algebra ${\cal A}_\ptl$ is usually referred to as a {\it Witt type Lie algebra} (e.g., \cite{SXZ,X}).

 For the Lie algebra
$\cal{A}_{\ptl}$, there are two natural ways to produce Lie coalgebras over
some subspaces of $\cal{A}^*,$ denoted as
$(\cal{A}_{\ptl})^{\circ}, (\cal{A^{\circ}})_{\ptl^{\circ}}$
respectively:\begin{enumerate}
\item The Lie coalgebra $(\cal{A}_{\ptl})^{\circ}$
is defined by \eqref{max-good} (cf.~Proposition \ref{PPPP}), with
cobracket, denoted $\D_\ptl$,  defined by \begin{equation*}
\D_{\ptl}(f) = (\mu (id \otimes \ptl - \ptl
\otimes id))^* (f)= (id \otimes \ptl - \ptl\otimes id)^*\mu^*(f)\mbox{ \ for $f \in (\cal{A}_{\ptl})^{\circ}$.}\end{equation*}
\item The Lie coalgebra
$(\cal{A^{\circ}})_{\ptl^{\circ}}$ is defined on the space ${\cal A}^\circ$, with cobracket, denoted $\D_{\ptl^\circ}$, being induced from the
coalgebra, i.e.
\begin{equation}\label{1111++}\D_{\ptl^{\circ}} (f) = (id \otimes \ptl^{\circ} -
\ptl^{\circ} \otimes id )\mu^{\circ} (f)\mbox{ \ for
$f \in \cal{A}^{\circ}$},\end{equation} where $\ptl^{\circ} =
\ptl^{\ast}|_{\cal{A^{\circ}}}$ is defined by $\ptl^{\circ} (f)(a) = f(\ptl(a))$ for $a \in \cal{A},
 $ and $\mu^{\circ} = \mu^{\ast}|_{\cal{A^{\circ}}}.$ 
\end{enumerate}
Note the difference that $(\cal{A}_{\ptl})^{\circ}$ is the maximal good subspace of the dual space $({\cal A}_\ptl)^*$ of the Lie algebra ${\cal A}_\ptl$, while
$(\cal{A^{\circ}})_{\ptl^{\circ}}$ as a space is the maximal good subspace of the dual space ${\cal A}^*$ of the associative algebra ${\cal A}$.

The following result (which holds for any field $\F$ of characteristic not $2$) can be found in \cite{6}.
\begin{prop}\adddot\label{theooo111}
Let $ \cal{A}$ be a commutative $\F$-algebra. 
Let $\ptl \in {\rm Der}_{\F}(\cal{A})$. Then
$(\cal{A}^{\circ})_{\ptl^{\circ}}$ is a Lie subcoalgebra of
$(\cal{A}_{\ptl})^{\circ}$. Furthermore, if there exist $c_1, c_2, \cdots, c_n
\in \cal{A}$ such that the ideal $I = (2\ptl(c_1), 2\ptl(c_2),
\cdots, 2\ptl(c_n) )$ has finite codimension, then
$(\cal{A}^{\circ})_{\ptl^{\circ}}=(\cal{A}_{\ptl})^{\circ}$.
\end{prop}

Now let $ \cal{A}= \F[x]$  ({\rm or} $\F[x^{\pm 1}]$) and  $\ptl=
\frac{d}{dx}$. Then we obtain the {\it one-side Witt algebra} ${\cal W}^+:=\F[x]_\ptl$
(or the {\it Witt algebra} $\cal{W}:=\F[x,x^{-1}]_\ptl$).
By Proposition \ref{theooo111},
$(\cal{A}^{\circ})_{\ptl^{\circ}}=(\cal{A}_{\ptl})^{\circ}$ in both cases.
The well-known Virasoro algebra $\cal{V}$ is the universal central extension of the
Witt algebra, namely, it has a basis $\{x^{m+1},c\,|\,m\in\Z\}$ with $c$ being a central element and
\begin{equation}\label{vir-rel}\mbox{$[x^{m+1}, x^{n+1}] =
(n\!-\!m)x^{m+n+1}\! +\! \frac{m^3 -m}{12}\d_{m+n,0}c\mbox{ \ for \ $m,n\in\Z$.}$}\end{equation}
It is proven in \cite{6} that the Virasoro coalgebra
$\cal{V^{\circ}}$ is isomorphic to $\cal{W^{\circ}}$, the dual of
Witt algebra. Thus from  the discussions above and Proposition \ref{theooo111}, we
obtain

\begin{theo}\adddot\label{theooo222}
Let $\cal{L}$ be the one side Witt algebra $\cal{W}^+$ $($resp., the Witt algebra $\cal{W}$
or the Virasoro algebra $\cal{V})$. Then the Lie coalgebra
$\cal{L^{\circ}}$ is equal to $(\F[x]^{\circ})_{\ptl^{\circ}}$ $($resp., $(\F[x,x^{-1}]^{\circ})_{\ptl^{\circ}})$, with
the~underlining~space~described~%
by~\eqref{poly}~$($resp.,~\eqref{lau-poly}$)$~and~%
the~cobracket~uniquely~determined~by
\begin{equation}
\label{dualLau-poly}
\D (\vp^n) = \mbox{$\sum\limits_{i+j=n+1}$} (j-i)\vp^i\otimes \vp^j\mbox{ \ \ for \ } n \in \N\ \mbox{$($resp.,
$n\in\Z)$},\end{equation}
where the sum is over $i,j\in\N$ $($resp., $i.j\in\Z{\sc\,}).$
\end{theo}

\begin{rema}\adddot\label{Fur-rrr}\rm\begin{enumerate}\item
Note that in case ${\cal L}={\cal W}$ or ${\cal V}$, the right-hand side of  \eqref{dualLau-poly} is always an infinite sum, which is an element in $({\cal L}\otimes{\cal L})^*$,
and its action on any element of $\F[x^{\pm 1}]
\otimes \F[x^{\pm 1}]$ can only produce finite many nonzero terms (cf.~\eqref{Dualll}).
Also note from Remark \ref{Poly-and-Laupoly} that \eqref{dualLau-poly} is not necessarily in ${\cal L}^\circ\otimes{\cal L}^\circ$ (which is the case when ${\cal L}={\cal W}$ or ${\cal V}$ because $\vp^n\notin\F[x,x^{-1}]^{\circ}$),
and that an element $f$ in $\F[x]^{\circ}$ or in $\F[x,x^{-1}]^{\circ}$ may be an infinite combination of $\vp^i$'s, thus $\D(f)$ is in general an infinite combination of elements of the
form 
\eqref{dualLau-poly} such that the resulting combination is in ${\cal L}^\circ\otimes{\cal L}^\circ$, which is guaranteed by Proposition \ref{PPPP}.
\item
We also remark that
in general, if elements $f=\sum_i f_i,\,g=\sum_jg_j$ are infinite sums in a Lie bialgebra,
it does not mean that \begin{equation}\label{Holdsss}\ \ \ \ \ \ \ \ \ \ \ \ \
\mbox{$\D(f)=\sum\limits_i\D(f_i)$  or $[f,g]=\sum\limits_{i,j}[f_i,g_j]$
if $f=\sum\limits_i f_i,\,g=\sum\limits_jg_j$ are infinite sums.}\end{equation} However in our case, since
we consider dual structures, elements $f,g$ are in certain dual spaces, because of the pairing like in \eqref{Dualll}, we see that \eqref{Holdsss} holds as long as it makes sense. Thus  in
 Theorem \ref{theooo222}, the cobracket is  uniquely determined by 
 \eqref{dualLau-poly}.
\end{enumerate}
%
\end{rema}

\ni{\it Proof of Theorem \ref{theooo222}.~}~We only need to prove the
cobracket relations.
We give a proof only for the case ${\cal L}={\cal W}^+$ as the proof for other cases is analogous. Thus let  $\vp^n \in $ $(\cal{W}^+)^{\circ}$ with $n\in\N$, and suppose
$\mu^{\circ} (\vp^n) =
\sum_{s,t\in\N} c_{st} \vp^s \otimes \vp^t$ for some $c_{st}\in\F.$ Then for $i,j\in\N$,
$$\mbox{$c_{ij}=
 \sum\limits_{s, t\in\N} c_{st}
\vp^s (x^i)\vp^t(x^j)=
\mu^{\circ}(\vp^n)(x^i \otimes x^j)=\vp^n(\mu(x^i \otimes x^j))=\vp^n(x^{i+j})=\d_{n,i+j}.
$}$$
Thus $\mu^{\circ}(\vp^n)\! =\! \sum_{i,j\in\N:\,i+j =n} \vp^i \otimes \vp^j$.
Similarly, if we assume $\ptl^{\circ}(\vp^n)\!=\!\sum_{s\in\N}c_s\vp^s$, then $c_i\!=\!
\ptl^{\circ}(\vp^n)(x^i)\!=\!\vp^n(\ptl(x^i))\!=\!i\d_{n,i-1}$, i.e., $\ptl^{\circ}(\vp^n)\!=\!(n\!+\!1)\vp^{n+1}$.
By Proposition \ref{theooo111} and \eqref{1111++}, we
\vspace*{-1pt}have
$$\mbox{$\D (\vp^n)\!=\! (id \otimes \ptl^{\circ}\! -\! \ptl^{\circ} \otimes id) \mu^{\circ}
 (\vp^n)\!=\!\sum\limits_{i+j=n} (id \otimes \ptl^{\circ} \!-\! \ptl^{\circ} \otimes id)(\vp^i \otimes
 \vp^j)\!=\!\sum\limits_{i+j=n+1} (j\!-\!i)\vp^{i}\otimes \vp^j,$}$$
\vspace*{-15pt}\

\noindent where the sums are over $i,j\in\N$, proving \eqref{dualLau-poly}.
\QED\vskip5pt
%
From now on,  $\cal{L}$ always denotes one of the one sided Witt algebra ${\cal W}^+$, the Witt algebra ${\cal W}$ and the
Virasoro algebra $\cal{V}$.
We summarize some results of \cite{9} as follows.
\begin{prop}\adddot\label{prppp00}
\begin{enumerate}
\item
Every Lie
bialgebra structure on ${\cal L}$
 is  coboundary triangular associated to a solution $r$ of CYBE \eqref{CYBE} of the form $r= a\otimes b- b\otimes a$ for some nonzero $a, b \in\cal{ L}$ satisfying $[a,b]=kb$ for some $0\ne k\in\F$.
\item
Let 
 $\gg $ be an infinite dimensional Lie
subalgebra of $\cal{W}$ such that $x \in \gg $ and $\gg  \ncong \cal{W}$
as Lie algebras. Let $\gg ^{(n)}$ be the Lie bialgebra structure on $\gg $
associated to the solution $r_n = x \otimes x^n - x^n \otimes x$
of the CYBE for any $x^n \in \gg $. Then every Lie bialgebra structure on
$\gg $ is isomorphic to $\gg ^{(n)}$  for some $n$ with $x^n \in \gg $.
\end{enumerate}\end{prop}


Using Theorem \ref{theooo222} and Proposition \ref{prppp00}, we can obtain the dual Lie bialgebra structures on ${\cal L}$ as follows.

\begin{theo}\adddot\label{dual-L}
%
%
Let
$(\cal{L}, [\cdot, \cdot], \D)$ be the coboundary triangular Lie bialgebra with
cobracket associated to the solution
$ r= x\otimes x^n -
x^n \otimes x$
of CYBE for some $1\ne n\in\Z$. Then the dual Lie
bialgebra
is $(\cal{L}^{\circ}, [\cdot, \cdot], \D)$
with the underlining space $\cal{L}^\circ$ and the  cobracket $\D$ determined by
Theorem $\ref{theooo222}$, and the bracket uniquely
determined by the following $($cf.~\eqref{Holdsss}$)$.\begin{enumerate}\item If $\cal{L}={\cal W}^+$, then $n\in\Z_+$ and for $i,j\in\Z_+$\vspace*{-5pt},
\begin{equation}\label{one-si}
[\vp^i,\vp^j]=\left\{\begin{array}{lll}
(2n\!-\!j\!-\!1)\vp^{j+1-n}\!\!&\mbox{if \ }i\!=\!1,\ 1\!\ne\!j\!\geq\! n\!-\!1,\\ 
(j-1)\vp^j&\mbox{if \ }i\!=\!n,\ j\!\neq \!1, i,n,\\ 
0&\mbox{if \ $i,j\notin\{1,n\}$ \ or \ $i\!=\!1,\ j\!<\!n\!-\!1$}.
\end{array}\right.
\end{equation}
\item If $\cal{L}={\cal W}$ or $\cal V$, then for $i,j\in\Z$\vspace*{-5pt},
\begin{equation}\label{twoone-si}
[\vp^i,\vp^j]=\left\{\begin{array}{lll}
(2n\!-\!j\!-\!1)\vp^{j+1-n}\!\!&\mbox{if \ }i\!=\!1,\ j\!\ne\! 1,\\ 
(j-1)\vp^j&\mbox{if \ }i\!=\!n,\ j\!\neq \!1,i, n,\\ 
0&\mbox{if \ $i,j\notin\{1,n\}$
}.
\end{array}\right.
\end{equation}
%
%
%
%
\end{enumerate}
\end{theo}

\ni{\it Proof.~}~By Theorem \ref{theooo222}, it remains to prove \eqref{one-si} and \eqref{twoone-si}. We only prove
\eqref{one-si} as the proof of \eqref{twoone-si} is similar.
Using notation
as in \eqref{Dualll}, by Proposition \ref{PPPP}, we obtain
\begin{eqnarray*}
 \langle [\vp^i, \vp^j], x^{m}\rangle &\!\!\!\!=\!\!\!\!& \langle \vp^i \otimes
\vp^j,
\D(x^{m})\rangle
=
\langle \vp^i \otimes \vp^j,x\cdot r \rangle \nonumber\\
&\!\!\!\!=\!\!\!\!& \langle \vp^i \otimes \vp^j, [x^{m}, x] \otimes x^{n} + x \otimes
[x^{m}, x^{n}]
- [x^{m}, x^{n}] \otimes x - x^{n} \otimes [x^{m},x]\rangle \nonumber\\
&\!\!\!\!=\!\!\!\!&(1-m)(\d_{i,m}\d_{j,n}-\d_{i,n}\d_{j,m})+(n-m)(\d_{i,1}\d_{j,m+n-1}-\d_{i,m+n-1}\d_{j,1})\nonumber\\
&\!\!\!\!=\!\!\!\!&\langle(1\!-\!i)\d_{j,n}\vp^i\!-\!(1\!-\!j)\d_{i,n}\vp^j\!+\!
(2n\!-\!j\!-\!1)\d_{i,1}\vp^{j+1-n}\!-\!(2n\!-\!i\!-\!1)\d_{j,1}\vp^{i+1-n},x^m\rangle.
\end{eqnarray*}
Thus
\begin{eqnarray}\label{lll}
 [\vp^i, \vp^j]=(1\!-\!i)\d_{j,n}\vp^i\!-\!(1\!-\!j)\d_{i,n}\vp^j\!+\!
(2n\!-\!j\!-\!1)\d_{i,1}\vp^{j+1-n}\!-\!(2n\!-\!i\!-\!1)\d_{j,1}\vp^{i+1-n}.
\end{eqnarray}
From this one immediately obtains $[\vp^i, \vp^j]=0$ if $i,j\in\{1,n\}$.
Now suppose $i=1\ne j$. Noting that $n\ne1$, we obtain from \eqref{lll} that $ [\vp, \vp^j]=
(2n\!-\!j\!-\!1)\vp^{j+1-n}$. In particular $ [\vp, \vp^j]=0$ if $j<n-1$ as in this case $\vp^{j+1-n}$ is not defined in $\F[x]^\circ$ (cf.~Proposition \ref{F[x]}). Thus we have the first and third cases of \eqref{one-si}.

Finally suppose $i=n$ and $j\ne 1,i,n$. Then \eqref{lll} gives $ [\vp^n, \vp^j]=
(j-1)\d_{i,n}\vp^j$, which is the second case of \eqref{one-si}. This completes the proof.
\QED
\begin{rema}\rm
Note from Proposition \ref{prppp00}$\sc\,$(2) that every nonzero solution $r$ of CYBE of ${\cal W}^+$ has the form $r=x\otimes x^n-x^n\otimes x$ for some $1\ne n\in\Z_+$.
Thus Theorem \ref{dual-L} in particular gives a complete classification of dual Lie bialgebras of ${\cal W}^+$.
However for the case $\cal{L}={\cal W}$ or $\cal V$,
a classification of
noncommutative $2$-dimensional Lie subalgebras of $\cal L$ (namely a classification of pairs $(a,b)$ of nonzero elements $a,b\in{\cal L}$ satisfying $[a,b]=kb$ for some $0\ne k\in\F$) is still an open problem, thus a classification of dual Lie bialgebras of $\cal L$ remains open.
\end{rema}

A series of $2$-dimensional Lie
subalgebras of $\cal{L}$ for ${\cal L}={\cal W}$ or $\cal V$ were presented in
\cite{11}. We present another series of $2$-dimensional Lie
subalgebras of $\cal{L}$ as follows.

\begin{lemm}\adddot\label{Lemmm}Let ${\cal L}=\cal{W}$ or $\cal V$, and fix $n\in\Z,\,\ell ,k\in\F$ with $n\ne1$ and $\ell ,k\ne0$. \vspace*{-5pt}Denote $$\mbox{$X = -\frac{1}{n-1}x+ \ell  x^n,\ \ \ \, Y
=-\frac{k}{2(n-1)\ell } x^{-n+2} + kx - \frac{(n-1)\ell k}{2}x^n \ \in\
\cal{L}\vspace*{-5pt}.$}$$ Then ${\rm span}_\F\{X, Y \}$ is a $2$-dimensional Lie subalgebra of $\cal{L}$,
and $[X, Y] = Y.$
\end{lemm}

\ni{\it Proof.~}~The proof is straightforward.\QED\vskip5pt
Now we are able to obtain the following main result of this section.
\begin{theo}\adddot\label{Last-Theo} Let $\cal{L}$, $X, Y$ be as in Lemma $\ref{Lemmm}$. Let
$(\cal{L},
$$[\cdot, \cdot], \D)$ be the coboundary triangular
Lie bialgebra associated to the solution
 $ r = X\otimes Y - Y \otimes X$ of CYBE.
Then the dual Lie bialgebra is $(\cal{L}^{\circ}, [\cdot, \cdot], \D)$
with the underlining space $\cal{L}^\circ$ and the  cobracket $\D$ determined by
Theorem $\ref{theooo222}$, and the bracket uniquely
determined \vspace*{-3pt}by  $($cf.~\eqref{Holdsss}$)$
%
\begin{equation}\label{one-si++++}
[\vp^i,{\ssc\!}\vp^j]\!=\!\!\left\{\!\!\!{\ssc\!}\begin{array}{lll} -k({\ssc\!}\d_{j, n}\!{\sc\!}+\!\d_{j,2-n}{\ssc\!})\vp^j
{\sc\!} \!-\!\frac{k(j+2n-3)}{2(n-1)^2\ell
}\vp^{j+n-1}{\sc\!}\!-\!
\frac{k\ell }{2}({\ssc\!}2n\!-\!j\!-\!1{\ssc\!})\vp^{j-n+1}\!\!\!&\mbox{if }i\!=\!1\!\ne\!j,\\[4pt]
-\frac{k}{2}\d_{j, n}\vp\!{\sc\!} -\!
\frac{k\ell }{2}({\ssc\!}n\!-\!1{\ssc\!})\d_{j, n}\vp^{2-n}{\sc\!} \!+\!\frac{k(j-1)}{2(n-1)^2\ell }\vp^j
{\sc\!}\!+\!\frac{k(2n-j-1)}{2(n-1)}\vp^{j-(n-1)}\!\!\!&\mbox{if }i\!=\!2\!-\!n,{\sc\,} j\!\ne\!1,2\!-\!n,\\[4pt]
\frac{k(j+2n-3)}{2(n-1)}\vp^{ j+n-1}
+\frac{k\ell }{2}(1-j)\vp^j&\mbox{if }i\!=\!n,{\sc\,} j\!\ne\!1,2\!-\!n,n,\\[4pt]
0&\mbox{if $i,j\notin\{1,2\!-\!n,n\}$}.
\end{array}\right.\!\!\!\!\!\!
\end{equation}
%
%
\end{theo}

\ni{\it Proof.~}~For convenience, we denote $\ell_0 \!=\!-\frac{1}{n-1}$,
 $k_0\! =\! -\frac{k}{2(n-1)\ell },$ $k_1\!=\!- \frac{(n-1)\ell k}{2}$. Thus
$X\!=\!\ell_0 x \!+\! \ell x^n$ 
and $Y = k_0x^{-n+2} + kx + k_1x^n $. \vspace*{-3pt}Then
\begin{eqnarray*}
\!\!\!\!\!\!\!\!\langle [\vp^i, \vp^j], x^m\rangle  &\!\!\!\!=\!\!\!\!& \langle \vp^i \otimes
\vp^j,
\D(x^m)\rangle
\\&\!\!\!\!=\!\!\!\!&
\langle\vp^i,[x^m, X\vp^j(Y)\!- \!Y\vp^j(X)]\rangle
+ \langle\vp^j,[x^m, Y\vp^i(X)\! -\! X\vp^i(Y)]\rangle\nonumber\\
&\!\!\!\!=\!\!\!\!&A_{ij}-A_{ji},
\end{eqnarray*}
\vspace*{-23pt}\

\noindent where\vspace*{-7pt},
\begin{eqnarray}\label{MMMM}
A_{ij}&\!\!\!\!=\!\!\!\!& \Big\langle\vp^i,[x^m, (\ell_0 x \!+\! \ell x^n) \vp^j(k_0 x^{2-n} \!+\! kx\! +\! k_1x^n)
\!-\! (k_0 x^{2-n}\! +\! kx\! + \!k_1x^n)\vp^j(\ell_0 x \!+\! \ell x^n)]\Big\rangle\nonumber\\
&\!\!\!\!=\!\!\!\!&\Big\langle
(k_0\d_{j,2-n}+k\d_{j,1}+k_1\d_{j,n})(\ell_0(1-i)\vp^i+\ell(2n-i-1)\vp^{i-n+1})\nonumber\\
&\!\!\!\!\!\!\!\!&\!-\!(\ell_0\d_{j,1}\!+\!\ell\d_{j,n})(k_0(3\!-\!i\!-\!2n)
\vp^{i+n-1}\!\!+\!k(1\!-\!i)\vp^i\!{\sc\!}+\!k_1(2n\!-\!i\!-\!1)\vp^{i-n+1}),x^m\Big\rangle.
\end{eqnarray}
From this, we see that $A_{ij}=0$ if $j\ne1,2-n,n$. Thus $[\vp^i,\vp^j]=0$ if $i,j\notin\{1,2-n,n\}$, which proves the last case of \eqref{one-si++++}.

Now assume $i=1\ne j$. By \eqref{MMMM}, we obtain
\begin{eqnarray*}
[\vp,\vp^j] &\!\!\!\!=\!\!\!\!&
(k_0\d_{j,2-n}+k_1\d_{j,n})\ell(2n-2)\vp^{2-n}-\ell\d_{j,n}(k_0(2-2n)\vp^{n}+k_1(2n-2)\vp^{2-n})
\\
&\!\!\!\!\!\!\!\!&-
k(\ell_0(1\!-\!j)\vp^j\!+\!\ell(2n\!-\!j\!-\!1)\vp^{j-n+1})
\\
&\!\!\!\!\!\!\!\!&+\ell_0(k_0(3\!-\!j\!-\!2n)\vp^{j+n-1}\!+
\!k(1\!-\!j)\vp^j\!+\!k_1(2n\!-\!j\!-\!1)\vp^{j-n+1}) \\
&\!\!\!\!=\!\!\!\!&\mbox{$-k\d_{j, 2-n}\vp^{2-n} \!-\! k \d_{j, n}\vp^n
\!-\!\frac{k\ell}{2}(2n\!-\!j\!-\!1)\vp^{j-n+1}
\!+\!\frac{k}{2(n-1)^2\ell}$}(3\!-\!j\!-\!2n)\vp^{j+n-1},
\end{eqnarray*}
which proves the first case of \eqref{one-si++++}.

Next assume $i=2-n$ and $j\ne1,2-n$. By \eqref{MMMM}, we obtain
\begin{eqnarray*}
[\vp^{2-n},\vp^j]&\!\!\!\!=\!\!\!\!&
k_1\d_{j,n}(\ell_0(n\!-\!1)\vp^{2-n}\!+\!\ell(3n\!-\!3)\vp^{3-2n})
\nonumber\\
&\!\!\!\!\!\!\!\!&
-\ell\d_{j,n}(k_0(1\!-\!n)\vp\!+\!k(n\!-\!1)\vp^{2-n}\!+\!k_1(3n-3)\vp^{3-2n})
\nonumber\\
&\!\!\!\!\!\!\!\!&- k_0(\ell_0(1-j)\vp^j+\ell(2n-j-1)\vp^{j-n+1})\\
&\!\!\!\!=\!\!\!\!&\mbox{$-\frac{k\ell }{2}(n-1)\d_{j, n}\vp^{2-n}
-\frac{k}{2}\d_{j,n}\vp +\frac{k(j-1)}{2(n-1)^2\ell}\vp^j +
\frac{k(2n-j-1)}{2(n-1)}\vp^{j-n+1}$},\end{eqnarray*}
which proves the second case of \eqref{one-si++++}.

Finally assume $i=n$ and $j\ne1,2-n,n$. By \eqref{MMMM}, we \vspace*{-5pt}obtain
\begin{eqnarray*}
[\vp^n,\vp^j]&\!\!\!\!=\!\!\!\!&
-k_1(\ell_0(1\!-\!j)\vp^j\!+\!\ell(2n\!-\!j\!-\!1)\vp^{j-n+1})
\nonumber\\
&\!\!\!\!\!\!\!\!&+
\ell(k_0(3\!-\!j\!-\!2n)\vp^{j+n-1}\!+\!k(1\!-\!j)\vp^j\!+\!k_1(2n-j-1)\vp^{j-n+1})\\
&\!\!\!\!=\!\!\!\!&\mbox{$\frac{k(2n+j-3)}{2(n-1)}\vp^{j+n-1} + \frac{k\ell}{2}(1-j)\vp^j ,$}
\end{eqnarray*}
which proves the third case of \eqref{one-si++++}.
This completes the proof of the theorem. \QED
\begin{rema}\rm\label{final-rem} Theorem \ref{Last-Theo} in particular gives a series of infinite dimensional Lie algebras (which seem to be new to us) on the space ${\mathcal S}:=\F[\vp,\vp^{-1}]$ with bracket defined by \eqref{one-si++++}, for various
 $n\in\Z,\,\ell ,k\in\F$ with $n\ne1$ and $\ell ,k\ne0$  (note from Remark \ref{Poly-and-Laupoly} that in fact ${\cal S}\cap {\cal L}^\circ=0$, however the bracket in \eqref{one-si++++} is well-defined on $\cal S$).
\end{rema}
We close the paper by proposing the following question: In which case will the dual Lie bialgebra of a coboundary triangular
Lie bialgebra be  coboundary triangular?

\small


\begin{thebibliography}{DWH99}
\def\re#1{\label{#1}\bibitem{#1}}
\small
\parindent=8ex\parskip=1.3pt\baselineskip=1.3pt




\re{1} R.E. Block,  Commutative Hopf algebras, Lie coalgebras, and divided powers, J. Algebra  96 (1985), 275--306.

\re{1-2} R.E. Block, P. Leroux, Generalized dual coalgebras of algebras, with applications to cofree coalgebras, J. Pure Appl. Algebra 36 (1985), 15--21.


\re{D}B. Diarra,
On the definition of the dual Lie coalgebra of a Lie algebra, Publ. Mat. 39 (1995), 349--354.


\re{2} V. Drinfel'd, Quantum groups, Proceedings ICM (Berkeley 1986),
Providence: Amer Math Soc, 1987, 789-820.

\re{G}G. Griffing,  The dual coalgebra of certain infinite-dimensional Lie algebras, Comm. Algebra 30 (2002), 5715--5724.


\re{3} S. Majid, Foundations of quantum group theory, Cambridge
University Press 1995.

\re{4} W. Michaelis, A class of infinite dimensional Lie bialgebras
containing the Virasoro algebras, Adv. Math. 107 (1994), 365--392.

\re{M1} W. Michaelis,
The dual Lie bialgebra of a Lie bialgebra, Modular interfaces (Riverside, CA, 1995), 81--93,
AMS/IP Stud. Adv. Math., 4, Amer. Math. Soc., Providence, RI, 1997.

\re{5}W.D. Nichols, The structure of the dual Lie coalgebra of
the Witt algebras, J. Pure Appl. Algebra 68 (1990), 359--364.

\re{6} W.D. Nichols, On Lie and associative duals, J. Pure Appl.
Algebra 87 (1993), 313--320.

\re{7} B. Peterson, E.J. Taft, The Hopf algebra of linearly
recursive sequences, Aequationes Mathematicae 20 (1980), 1--17,
University Waterloo.

\re{8} E.J. Taft, Witt and Virasoro algebras as Lie bialgebras,
 J. Pure Appl. Algebra 87 (1993), 301--312.


\re{9} S.-H. Ng, E.J. Taft, Classification of the Lie bialgebra
structures on the Witt and Virasoro algebras, J. Pure Appl. Algebra
151 (2000), 67--88.

\re{10} G. Song, Y. Su, Lie Bialgebras of generalized Witt type,
Science in China: Series A Mathematics  49(4) (2006), 533--544.

\re{SXZ}Y. Su, X. Xu, H. Zhang,
Derivation-simple algebras and the structures of Lie algebras of Witt type, J. Algebra 233 (2000), 642--662.

\re{11} Y. Su, K. Zhao, Generalized Virasoro and super-Virasoro
algebras and modules of the intermediate series, J. Algebra 252
(2002), 1--19.

\re{12} M.E. Sweedler, Hopf Algebras, W. A. Benjamin, Inc. New
York, 1969.

\re{X} X. Xu, New generalized simple Lie algebras of Cartan type over a field with characteristic
0, J. Algebra 224 (2000), 23--58.

\re{13} Y. Wu, G. Song, Y. Su, Lie bialgebras of generalized
Virasoro-like type, Acta Math. Sinica Engl. Ser. 22 (2006),
1915--1922.

\re{14}B. Xin, G. Song, Y. Su, Hamiltonian type Lie bialgebras,
Science in China Series A: Mathematics 50 (2007),
1267--1279.

\end{thebibliography}
\end{document}